%% file: Fibers1.tex
\documentstyle[twoside,12pt,toc_small,bib_toc,psfig]{article}

\input{LaTeX-top}

\def\M{{\bf M}}

\def\LabelCaption#1#2{
   \label{#1}
   \centerline{\parbox{\captionwidth}{   %\vspace{-10mm} 
   \caption{\sl #2}  } } }
\newfont{\bbflarge}{msbm10 at 21pt}
\def\Clarge{\mbox{\bbflarge C}}

\edef\theta{\vartheta}
\long\def\longhide#1{}

\def\MiLemRiesz{15.3}  % The Riesz Theorem in MiIntro
\def\MiLemRieszA{A.3}  % Nochmal
\def\MiCor{18.6}       % Kein Strahl bei Cremer Punkt
\def\MandelBranch{2.2} % The Branch Theorem
\def\LemApproxPeriodic{3.5} % Approximation of periodic ray pairs
\def\LemApproxPreperiodic{3.9} % same, preperiodic case
\def\ThmRayCorrespondence{2.1} % Ray transfer dyn - para
\def\CorThreeRaysOneFiber{3.16} % Three rays in para space

\def\newsection#1{
   \section{#1}
}

\newlength\captionwidth
\title{\vspace{-15mm} On Fibers and Local Connectivity\\
of Compact Sets in $\Clarge$}
\author{Dierk Schleicher\\
Technische Universit\"at M\"unchen}
\date{}

\begin{document}

\maketitle
\thispagestyle{empty}
\input{imsmark.tex}
\def\IMSmarkvadjust{-30pt}
\SBIMSMark{1998/12}{December 1998}{}

%\centerline{\reminder{Draft of \today} }

\begin{center}
\begin{minipage}{110mm}
\def\toc_vspace{1.5em}
\def\tocname{Contents}
\tableofcontents{0.3em}
\end{minipage}
\end{center}
\vspace{5mm}

\begin{abstract}  
A frequent problem in holomorphic dynamics is to prove local
connectivity of Julia sets and of many points of the Mandelbrot set;
local connectivity has many interesting implications. The intention of
this paper is to present a new point of view for this problem:
we introduce {\em fibers} of these sets (Definition~\ref{DefFiber}),
and the goal becomes to show that fibers are ``trivial'', i.e.\ they
consist of single points. The idea is to show ``shrinking of puzzle
pieces'' without using specific puzzles. This implies local
connectivity at these points, but triviality of fibers is a somewhat
stronger property than local connectivity. Local connectivity proofs
in holomorphic dynamics often actually yield that fibers are trivial,
and this extra knowledge is sometimes useful. 

Since we believe that fibers may be useful in further situations,
we discuss their properties for arbitrary compact connected and
full sets in the complex plane. This allows to use them for connected
filled-in Julia sets of polynomials, and we deduce for example that
infinitely renormalizable polynomials of the form $z^d+c$ have the
property that the impression of any dynamic ray at a rational angle
is a single point. An appendix reviews known topological properties
of compact, connected and full sets in the plane. 

The definition of fibers grew out of a new brief proof that the
Mandelbrot set is locally connected at every Misiurewicz point and at
every point on the boundary of a hyperbolic component. This proof
works also for ``Multibrot sets'', which are the higher degree cousins
of the Mandelbrot set. These sets are discussed in a self-contained
sequel \cite{FiberMandel}. Finally, we relate triviality of fibers to
tuning and renormalization in \cite{FiberTuning}.

\end{abstract}

\newpage

\newsection {Introduction}
\label{SecIntro}

A great deal of work in holomorphic dynamics has been done in recent
years trying to prove local connectivity of Julia sets and of many
points of the Mandelbrot set, notably by Yoccoz, Lyubich, Levin,
van Strien, Petersen and others. One reason for this work is that the
topology of Julia sets and the Mandelbrot set is completely described
once local connectivity is known. Another reason is that local
connectivity of the Mandelbrot set implies that hyperbolicity is
dense in the space of quadratic polynomials, and that the dynamics can
completely be classified by its combinatorics plus multipliers of
attracting orbits.

In this paper, we introduce {\em fibers} of Mandelbrot and Julia sets
and  shift the focus from local connectivity to a closely related but
somewhat stronger concept which we call triviality of fibers. It can be
observed that people often prove that fibers are trivial when they only
speak about local connectivity. However, triviality of fibers has quite
a few useful properties: it allows to draw some conclusions which do
not follow from local connectivity, and it makes several proofs more
transparent. On the other hand, the concept of trivial fibers is not
too restrictive: every compact connected and full subset of the complex
plane which is locally connected has only trivial fibers for an
appropriate choice of external rays used in the construction of fibers. 

A fundamental construction in holomorphic dynamics is called the {\em
puzzle}, introduced by Branner, Hubbard and Yoccoz. A typical proof
of local connectivity consists in establishing {\em shrinking of
puzzle pieces} around certain points. This is exactly the model for
fibers: the fiber of a point is the collection of all points which
will always be in the same puzzle piece, no matter how the puzzle was
constructed. Our arguments will thus never use specific puzzles. 

This paper is the first in a series: we introduce fibers and discuss
their properties for arbitrary compact connected and full subsets of
the complex plane. In particular, we explain the relation between
triviality of fibers, local connectivity and landing properties of
external rays (Section~\ref{SecFibers}). It turns out that it is
possible to construct certain bad subsets of $\C$ for which fibers
behave rather badly. However, we will give a criterion in
Lemma~\ref{LemFibersNice} which will ensure that fibers are
well-behaved, and this criterion will usually be satisfied in
holomorphic dynamics.

In Section~\ref{SecJulia}, we apply fibers to connected filled-in
Julia sets of polynomials and show that they are generally quite
well-behaved. As a new result, we show that many Julia sets have the
property that all periodic external rays have impressions consisting
only of their landing points. These Julia sets include infinitely
renormalizable Julia sets of polynomials with a single critical point. 
We will need Thurston's No Wandering Triangles Theorem, which we cite
here with his proof and permission. 

The paper concludes with an appendix about compact connected full
(and sometimes locally connected) subsets of the complex plane.
Several well known results which are needed elsewhere in the paper are
collected there, often with proofs included for easier reference.

In \cite{FiberMandel}, we will use fibers to give a new proof that the
Mandelbrot set (and more generally Multibrot sets) have trivial
fibers at Misiurewicz points and at all boundary points of hyperbolic
components, including roots of primitive components. An immediate
corollary will be local connectivity at these points. Finally in
\cite{FiberTuning}, we will discuss how triviality of fibers is
related to renormalization and tuning: in parameter space, it is
preserved under tuning, and any Julia set of the form $z\mapsto z^d+c$
has all its fibers trivial if and only if any of its renormalizations
has this property; again, the same follows for local connectivity of
these sets.

\longhide{ ******
We discuss polynomials of the form $z\mapsto z^d+c$, for
arbitrary complex constants $c$ and arbitrary degrees $d\geq 2$.
These are, up to normalization, exactly those polynomials which have
a single critical point. Following a suggestion of Milnor, we call
these polynomials {\em unicritical} (or unisingular). We will always
assume unicritical polynomials to be normalized as above, and the
variable $d$ will always denote the degree. We define the {\em Multibrot
set} of degree $d$ as the connectedness locus of these families, that is
\[
\M_d := \{c\in\C: \mbox{ the Julia set of $z\mapsto z^d+c$ is
connected }\} \,\,.
\]
In the special case $d=2$, we obtain quadratic polynomials, and
$\M_2$ is the familiar {\em Mandelbrot set}. All the Multibrot sets
are symmetric with respect to the real axis, and they have also
$d-1$-fold rotation symmetries (see \cite{Intelligencer} with pictures
of several of these sets). The present paper can be read with the
quadratic case in mind throughout. However, we have chosen to do the
discussion for all the Multibrot sets because this requires only
occasional slight modifications -- and because recently interest in the
higher degree case has increased; see e.g.\ Levin and van
Strien~\cite{LvS}.

One problem in holomorphic dynamics is that many results are
folklore, with few accessible proofs published. In particular, many
of the fundamental results about the Mandelbrot set, due to Douady
and Hubbard, have been described in their famous ``Orsay notes''
\cite{Orsay} which have never been published. They are no longer
available, and it is not always easy to pinpoint a precise reference
even within these notes.

This paper provides proofs of certain key results about the
Mandelbrot sets, and more generally for Multibrot sets. Most of these
results are known at least for $d=2$, but several of our proofs are new
and sometimes more direct that known proofs. 

In a forthcoming paper \cite{FiberTuning}, we will give various
applications of the concept of fibers which are related to the concept of
{\em tuning}: triviality of fibers is preserved under tuning, both for
Julia sets and for Multibrot sets, and we will show that the Mandelbrot
set comes quite close to being arcwise connected. These statements, in
turn, have various further applications. 

The present paper is intended to serve several purposes: to give new
proofs for certain key results about the Mandelbrot set, to provide
precise references for the statements, to generalize them to the
Multibrot sets, and to identify and introduce fibers as a framework for
further work in a similar direction. Despite these rather different
purposes, we feel that the paper has an inner unity. Its organization is
as follows: in Section~\ref{SecFibers}, we introduce the concept of
fibers of arbitrary compact connected and full subsets of the complex
plane and show several fundamental results for these fibers. In
particular, we discuss the relation between triviality of fibers and
local connectivity.

The Mandelbrot and Multibrot sets are discussed in
Section~\ref{SecMultibrot}. We review certain fundamental results about
these sets which will later be needed, and we give some applications of
fibers. We show that the fiber of an interior point is trivial if and
only if it is in a hyperbolic component. Moreover, we show that local
connectivity is equivalent to triviality of all fibers (constructed using
rational parameter rays), and both conditions imply density of
hyperbolicity (using an argument of Douady and Hubbard). 

In Section~\ref{SecMandelFibers}, we prove that every Multibrot set has
trivial fibers and is thus locally connected at every boundary point
of a hyperbolic component and at every Misiurewicz point. This shows that
fibers of Multibrot sets have particularly convenient properties. We
also compare fibers to combinatorial classes.

In Section~\ref{SecJulia}, we apply the concept of fibers to
Julia sets. We will need Thurston's No Wandering Triangles Theorem,
which we cite here with Thurston's proof and permission. This allows to
conclude that for connected Julia sets of unicritical polynomials, local
connectivity is equivalent to all the fibers being trivial for a
specified collection of dynamic rays used in the construction of fibers.
A new result is that all periodic points of infinitely renormalizable
unicritical polynomials have trivial fibers, which implies that the
impression of any dynamic ray at a rational external angle is a single
point. This section is independent from the previous two sections.

The paper concludes with an appendix about compact connected full
and locally connected subsets of the complex plane. Several well
known results which are needed elsewhere in the paper are collected
there, often with proofs included for easier reference.

} %end of \longhide

{\sc Acknowledgements.}
This entire paper would not be without many interesting and useful
discussions with Misha Lyubich, which is just one reason why I am
most grateful to him. Thanks go also to Genadi Levin for an inspiring
discussion at an early stage of this paper. Most of this happened
during the special semester in 1995 at the Mathematical Sciences
Research Institute in Berkeley, which deserves special thanks for its
hospitality. Most of the writing of this paper was done in the
stimulating atmosphere at the Institute for Mathematical Sciences in
Stony Brook. I am most grateful to John Milnor for the invitation and
for many inspiring discussions and suggestions, including helpful
remarks after a first reading of the manuscript. Quite a few
interesting discussions with Saeed Zakeri and Adam Epstein were also
most helpful.

\newsection{Fibers and Local Connectivity} 
\label{SecFibers}

Our goal in this section is to introduce fibers of compact connected
and full subsets of $\C$. Fibers will be the topological building
blocks. We will discuss triviality of fibers and local connectivity
as two closely related concepts which are the focus of interest of a
lot of work, related for example to the Mandelbrot set.

Throughout this section, let $K$ be a connected, compact and full
subset of $\C$ (``full'' means that the complement has no bounded
components). External rays of $K$ are defined as inverse images of
radial rays under the Riemann map sending the exterior of $K$ to the
exterior of the unit disk, normalized so as to fix $\infty$ with
positive real derivative (in the special case that $K$ has conformal
radius one, this means that the Riemann map is tangent to the identity
at $\infty$). When dealing with dynamic and parameter planes, we will
speak of ``dynamic rays'' and ``parameter rays'' instead of external
rays.

\begin{definition}[Limit Set and Impression of External Ray]
\label{DefImpress} \lineclear 
We denote the external ray of $K$ at angle $\theta$ by $R_K(\theta)$. 
Its\/ {\em limit set} is
$L_K(\theta):=\ovl{R_K(\theta)}\cap K$. 
The {\em impression} of the ray is the set
\[ 
I_K(\theta):=
\bigcap_{\eps>0}\,\,\ovl{\bigcup_{|\phi-\theta|<\eps}L_K(\phi)}
\,\,. 
\]
We say that the external ray at angle $\theta$ {\em lands} if
its limit set is a single point.
\end{definition}
Equivalently, limit set and impression can be defined as the sets of all
possible limits 
\[
L_K(\theta) = \lim_{r\searrow 1} \Phi_K^{-1}(r e^{2\pi i\theta})
\,\,\mbox{ and }\,\,
I_K(\theta) = \lim_{r\searrow 1, \phi\to\theta}  
\Phi_K^{-1}(r e^{2\pi i\phi}) \,\,,
\]
where $\Phi_K\colon\Cbar-K\to \Cbar-\diskbar$ is the normalized Riemann
map. The impression obviously contains the limit set, and both are
compact, connected and non-empty. It may well happen that an external
ray lands while its impression is a continuum. As usual, we measure
external angles in full turns so that they live in $\Circle=\R/\Z$.
Let $Q\subset\Circle$ be any {\em countable} subset of angles such that
all the external rays at angles in $Q$ land. (One could allow larger
sets $Q$, for example the set of all angles such that the corresponding
rays land; by Fatou's Theorem, this set has full measure in $\Circle$.
However, in all the applications we have in mind, the set $Q$ will be
countable anyway, and the countability hypothesis makes a few arguments
more convenient; see also the remarks after Lemma~\ref{LemFibersNice}.)
In most cases, $Q$ will be the set of rational angles, in particular
when discussing Multibrot sets and {\em monic} polynomials (however,
when there are Siegel disks, we need to enlarge $Q$). We will often
loosely speak of an ``external ray in $Q$'' when we mean an external ray
such that its external angle is in $Q$, thus identifying rays with their
angles. 

The landing properties of external rays are studied by Carath\'eodory
Theory, which investigates into which pieces the boundary of $K$ can
be cut by external rays landing there (see for example
Milnor~\cite[Sections 15 and 16]{MiIntro}; recently, 
Petersen~\cite{LundePrimeEnds} has refined the study of this theory). We
are going to do a related study here, but we will look at the set $K$
from inside as well as from outside. 

\begin{definition}[Separation Line]
\label{DefSeparation} \lineclear
A {\em separation line} will be two external rays with angles in $Q$
which land at a common point on $\partial K$, or two such rays which
land at different points, together with a simple curve in the interior
of\/ $K$ connecting the two landing points. Two points $z,z'\in K$
{\em can be separated\/} if there is a separation line $\gamma$
avoiding $z$ and $z'$ such that these two points are in different
connected components of\/ $\C-\gamma$.
\end{definition}
The separation line should also contain the landing points of the
two rays. The curve in the interior of $K$ must land at the same
points as the two rays (where landing is understood in the same
sense as for rays). Therefore, any separation line will cut the
complex plane into two open parts. We will use these lines to
define fibers of $K$ and to construct connected neighborhoods of a
point when proving local connectivity at this point. 

When an interior component of $K$ is equipped with an arbitrary base
point, one might require the separation line within this component to be
the union of two ``internal rays'': since the interior component is
simply connected, there is a Riemann map from the component to $\disk$
sending the base point to $0$, and this map is unique up to rotation.
Internal rays are then inverse images of radial lines, and by
Lindel\"of's Theorem~\ref{ThmLindeloef} in the appendix, any point which
is accessible by a curve is in fact the landing point of a ray. As far
as the boundary of $K$ is concerned, there is nothing lost in restricting
to internal rays. We will not need to make this restriction. 

\begin{definition}[Fibers and Triviality]
\label{DefFiber} \lineclear
For any point $z\in K$, consider the set of points in $K$ which cannot
be separated from $z$. In this set, the connected component containing
$z$ will be called the {\em fiber} of $z$. We say that it is {\em
trivial} if it consists of the point $z$ alone. 
\end{definition}
\remark
Mandelbrot, Multibrot, and Julia sets are often studied with the help
of a partition called a ``Branner-Hubbard-Yoccoz jigsaw puzzle'', and
a lot of work is devoted to showing that ``puzzle pieces shrink to
points'' (expressed by Douady as ``points are points''). The idea
behind fibers is to capture the essential properties of these puzzles
without using any details about the exact construction of the puzzle. 
Douady's joke ``points are points'' can then be replaced by the more
precise (but dull) ``fibers are points''.

For the Mandelbrot and Multibrot sets, fibers and their triviality are
related to combinatorial classes and combinatorial rigidity. These
differ exactly at hyperbolic components: entire hyperbolic components
form combinatorial classes together with part of their boundaries,
while we want to distinguish all their points; compare
\cite{FiberMandel}. For parameter spaces, another way of saying
that the fiber of a point is trivial is that the space is ``fiber
rigid'' at this point. In order to avoid overusing the word
``rigidity'', we have decided not to use it for general sets $K$ or
for Julia sets and rather speak of ``triviality of fibers''.

\remark
For the given definition of fibers, it is possible to construct
compact connected full sets $K\subset\C$ for which fibers behave
badly. However, for the applications we have in mind, fibers usually
have quite nice properties because we can choose $Q$ so that the
landing points of the selected rays have trivial fibers: see
Lemma~\ref{LemFibersNice} and the remarks thereafter. We have not been
successful in finding a satisfactory definition of fibers which has
similar pleasant properties from the start for arbitrary sets $K$
without becoming too complicated for the sets we are interested in. 
Most of the problems are related to interior components of $K$. If
there is no interior, which is the case for many interesting Julia
sets, the situation generally becomes quite a bit easier.

In the following paragraph, we will describe some ``bad''
possibilities of fibers for appropriately constructed sets $K$, in
order to show what we need to have in mind in our proofs. 

For certain points $z$, it may happen that the set of points in $K$
which cannot be separated from $z$ is disconnected. This occurs whenever
there is a component $U$ of the interior of $K$ which has exactly two
boundary points $p,q$ each of which is accessible from within $U$
and the landing point of a ray in $Q$ (see Figure~\ref{FigFiber2Acc}).
Then every point in $U$ can be separated from any other point in $\ovl
U$ except $p$ and $q$. With our definition, the fiber of every point in
$U$ is trivial, while the fibers of $p$ and $q$ contain all of $\ovl U$.
This same example also shows that the relation ``$z_1$ is in the fiber
of $z_2$'' need not be symmetric or transitive. Fibers of different
points may also intersect without being equal: as an example, take the
filled-in Julia set of $z^2-1$ (the ``Basilica'') and let $Q$ be the
set of external angles of the form $a/(3\cdot 2^k)$ for integers
$k\geq 0$ and $a>0$; these are exactly the angles of external rays
landing together with another ray. Not allowing separation lines
through the interior of $K$ (or changing the topology so that curves
in the interior of $K$ cannot land at landing points of rays in $Q$),
then the fiber of any interior point in $K$ is the closure of its
connected component of the interior of $K$, and two such components
may have intersecting closures. --- Fibers may thus have some rather
unpleasant properties. One could try to remedy this by defining new
fibers to be the smallest equivalence classes which are topologically
closed and which contain entire fibers in the sense above. However, it
would then be possible that some point could be separated from every
other point without its fiber being trivial; in fact, the fiber could
be all of $K$. As mentioned above, in our applications we will usually
be able to choose our rays in $Q$ so that the fibers of their landing
points are trivial, and all these problems disappear by 
Lemma~\ref{LemFibersNice}. 

We begin by collecting a couple of useful properties of fibers which
are true in general.

\hide{
All these problems disappear when the landing points of the
rays in $Q$ have trivial fibers (see Lemma~\ref{LemFibersNice} and the
remarks thereafter). This ``nice'' situation is what we have in mind,
and it will usually be satisfied in our applications. We have not been
successful in finding a satisfactory definition of fibers which has
similar pleasant properties from the start for arbitrary sets $K$.
}

\begin{figure}[htbp]
%\vspace{68mm}
%\hskip 44mm 
%\special{picture Fiber2Acc scaled 500} \hfill\hfill
\centerline{
\psfig{figure=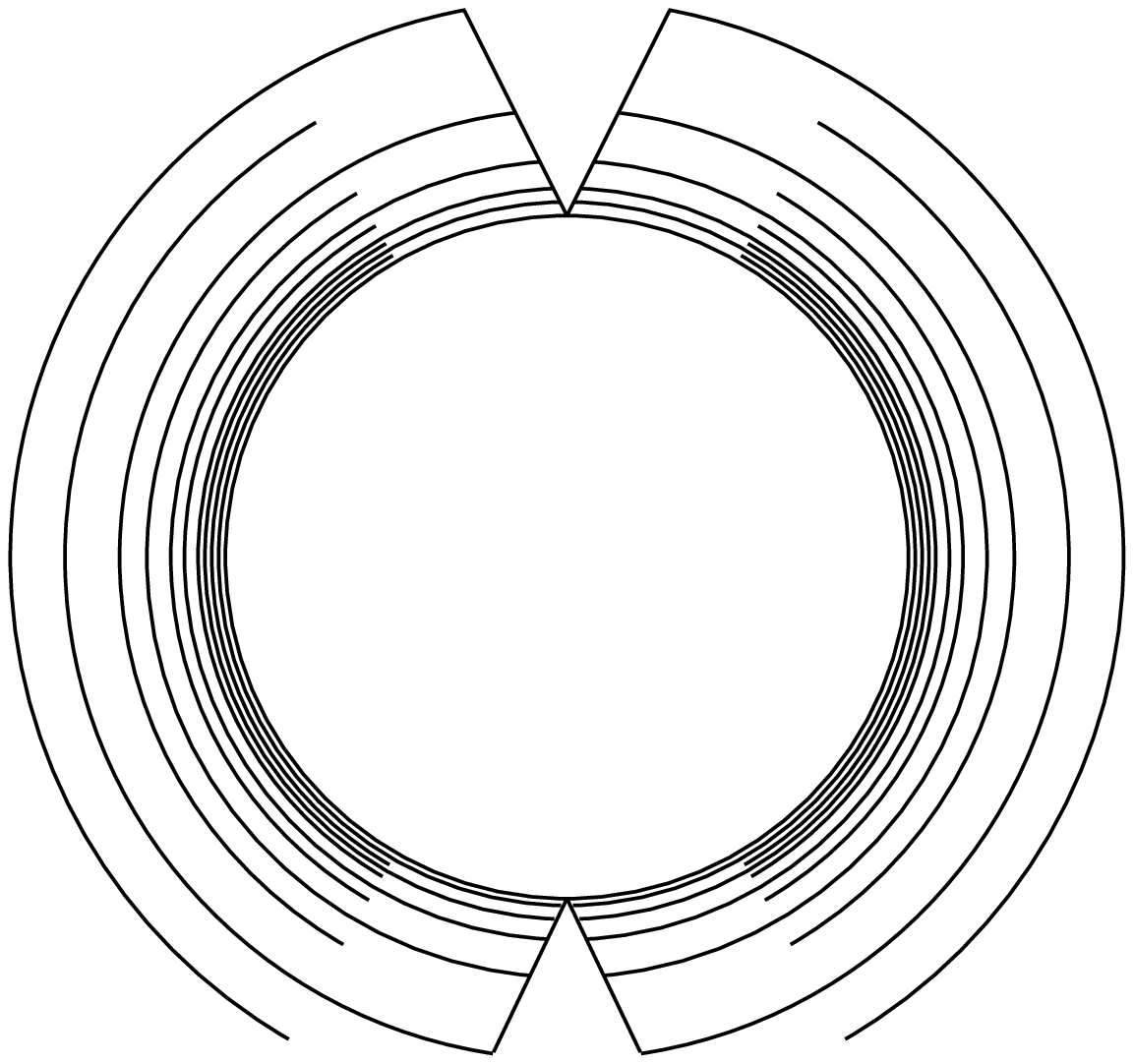,width=60mm}
}
\LabelCaption{FigFiber2Acc}
{A compact connected and full set $K$ such that an interior component
has exactly two boundary points which are simultaneously accessible
from inside and outside. Any separation line through the interior must
pass both of these boundary points, so these two boundary points
cannot be separated from any interior point.} \hfill
\end{figure}

\begin{lemma}[Properties of Fibers]
\label{LemFibers} \lineclear
Fibers have the following properties:
\begin{enumerate}
\item
Every fiber is compact, connected, and full.
\item
Any connected component of the interior of $K$ is either contained
in a single fiber, or the fiber of each of its points is trivial;
the latter happens if and only if at least two rays in $Q$ land on the
boundary of this connected component such that their landing points are
accessible from the inside of the component. 
\item
Let $z$, $z'$ be two points in $K$. If $z'$ is not in the fiber of $z$,
then $z$ and $z'$ can be separated using a separation line avoiding $z$
and $z'$, except in the following situation: $z$ is in the interior of
$K$, the interior component of $K$ containing $z$ has exactly two
boundary points which are landing points of rays in $Q$ and can be
connected by a curve in the interior of $K$, and both of these boundary
points have non-trivial fibers.
\item
Finally, if a fiber consists of more than a single point, then its
boundary is contained in the boundary of $K$.
\end{enumerate}
\end{lemma}
\proof
Any separation line avoiding $z$ obviously separates an open subset of
$K$ from $z$, so the fiber of $z$ is an intersection of closed sets and
thus closed. It is connected by definition. 

Let $U$ be a connected component of the interior of $K$. If at least
two rays with angles in $Q$ land on the boundary of $U$ such that there
is a curve in $U$ landing at the landing points of the rays, then every
point in $U$ can easily be separated from any other point in $U$, and
the fiber of every point in $U$ is trivial. Conversely, if not all of
$U$ is in the same fiber, there must be a curve in $U$ connecting two
landing points of rays in $Q$. But then the fiber of every point in $U$
is trivial. (However, fibers of boundary points of $U$ might not be
trivial.) 

Let $z$ be an arbitrary point in $K$ and assume first that the set of
separation lines is countable. Let $K_n$ be the closure of the connected
component of $K$ containing $z$ in the complement of the first $n$
separation lines avoiding $z$. It is compact and connected. It is also
full because the complements of $\C-K_n$ are open and connected, and
their unions are then also open and connected.  The fiber of $z$ is the
nested intersection of all the $K_n$ and thus full for the same reason. 

Any point $z'\in K-K_n$ can then be separated from $z$ be a separation
line avoiding both $z$ and $z'$, unless two separation lines meet more
than once and together separate a point from $z$ which is not separated
by any single separation line. But this can happen only if $z$ is in the
interior of $K$ and has a trivial fiber. Let $U$ be the connected
component of the interior of $K$ containing $z$. If at least three
boundary points of $U$ are landing points of rays in $Q$ and of curves
from within $U$, then $z$ can easily be separated from any point in
$K-\{z\}$ by a separation line avoiding both points. If this does not
happen, but the fiber of $z$ is trivial, there must be exactly two such
boundary points. But then the fibers of these boundary points must both
contain all of $U$. 

Since we have assumed the set $Q$ of rays to be countable, the choice of
the two external rays used for a separation line is also countable. Any
two separation lines using the same two external rays must either
coincide, or they must traverse the same interior component of $K$.
Therefore, a countable collection of separation lines is always
sufficient, and the proof above works in general. All fibers are thus
full. 

Now suppose that a point $z$ is a boundary point of a fiber $Y$ and
an interior point of $K$. Then the connected component of the
interior of $K$ containing $z$ must contain a non self-intersecting
curve connecting the landing points of two rays in $Q$, and $Y=\{z\}$.
\qed

\begin{lemma}[Impression is in Single Fiber] 
\label{LemImpressFiber} \lineclear 
For an external ray which lands (in particular for rays in $Q$),
the impression is con\-tai\-ned in the fiber of its landing point.
For a ray which is not in $Q$ (even if it does not land), the
impression is contained in the fiber of any point in the
impression.
\end{lemma}
\proof
For a ray in $Q$, let $z$ be its landing point; for a ray not in $Q$,
let $z$ be any point in the impression. Then $z\in\partial K$. We want
to show that any point $z'\in K$ which is not in the fiber of $z$ cannot
be in the impression of the ray. But this is obvious because $z$ and
$z'$ are separated by a separation line (Lemma~\ref{LemFibers}), and no
impression can extend over this separation. 
\qed
\remark
It is not quite true that the impression is contained in the fiber of
any point from the impression: if a ray in $Q$ is part of a
separating ray pair, then the impression may extend over both sides
of the separation line, while fibers of points from different sides
cannot contain each other. However, the fiber of the landing point
will still contain the entire impression.

\begin{lemma}[Boundary Points are in Impression]
\label{LemBdyInImpression} \lineclear
Every boundary point of $K$ is in the impression of at least one
external ray. If the fiber of a boundary point is trivial, then at
least one external ray lands there.
\end{lemma}
\proof
Let $z\in\partial K$ and let $(z_n)$ be a sequence of points in $\C-K$
tending to $z$. The external angles of $(z_n)$ must then have at
least one limit $\theta\in\Circle$, so that $z$ is in the impression
of the ray at angle $\theta$. If the fiber of $z$ is trivial, then
$z$ can be separated from any $z'\in K$ and the impression of the ray is
$\{z\}$, which implies in particular that the ray lands at $z$.
\qed

The next lemma shows that fibers behave particularly nicely if the
rays in $Q$ land at points with trivial fibers. 

\begin{lemma}[When Fibers Behave Nicely]
\label{LemFibersNice} \lineclear
If the landing points of all the rays in $Q$ have trivial fibers, then
the fibers of any two points are either equal or disjoint, and the
set $K$ splits into\/ {\em fibers} as equivalence classes of points
with coinciding fibers. In that case, there is an obvious map from
external angles to fibers of $K$ via impressions of external rays. This
map is surjective onto the set of fibers meeting $\partial K$
\end{lemma}
\proof
The relation ``$z_1$ is in the fiber of $z_2$'' is always reflexive.
When the landing points of rays in $Q$ have trivial fibers, then this
relation is also symmetric by Lemma~\ref{LemFibers}. In order to show
transitivity, assume that two points $z_1$ and $z_2$ are both in the
fiber of $z_0$. If they are not in the fibers of each other, then the
two points can be separated by a separation line avoiding $z_1$ and
$z_2$ (Lemma~\ref{LemFibers}). If such a separation line can avoid
$z_0$, then these two points cannot both be in the fiber of $z_0$. The
only separation between $z_1$ and $z_2$ therefore runs through the point
$z_0$, so $z_0$ cannot be in the interior of $K$ and rays in $Q$ land at
$z_0$. By assumption, the fiber of $z_0$ consists of $z_0$ alone. Any
two points with intersecting fibers thus have indeed equal fibers. The
map from external angles to fibers exists by
Lemma~\ref{LemImpressFiber}. It is surjective by
Lemma~\ref{LemBdyInImpression}. 
\qed
\remark
The situation described in this lemma is what we want fibers to be:
we want to speak of ``fibers of $K$'' rather than having to specify
which point of $K$ any fiber is seen from. This is one reason not to
make $Q$ unnecessarily large, or it would be harder to establish this
``nice'' property. We will show in \cite{FiberMandel} that
the Mandelbrot and Multibrot sets have this property: this amounts to
showing that they have trivial fibers  at the boundary of hyperbolic
components (including the roots of primitive components) and at
Misiurewicz points. Also, most Julia sets have ``nice'' fibers
(Section~\ref{SecJulia}). 

If the set $K$ has trivial fibers at all the landing points of rays in
$Q$, then it is not hard to show that the quotient of $K$ by
identifying points with coinciding fibers is a compact connected
locally connected Hausdorff space (for the proof of local
connectivity, see the proof of Proposition~\ref{PropFiberLocConn}
below). In fact, the topological pair $(\Cbar,K)$ modulo this
equivalence relation is homeomorphic to the topological pair
$(\Sphere,K')$ for a compact connected locally connected set $K'$:
this is due to Moore's Theorem assuring exactly that (see
Daverman~\cite{Dav}).

\begin{definition}[Local Connectivity] 
\label{DefLocConn} \lineclear 
A compact connected set $K\subset\C$ is called\/ {\em locally
connected at a point $z\in K$} if every neighborhood of $z$ contains
a subneighborhood intersecting $K$ in a connected set. If this
subneighborhood can always be chosen open, then\/ $K$ is said
to be {\em openly locally connected at\/ $z$}. We say that {\em $K$
is locally connected} if it is locally connected at each of its
points.
\end{definition}
\remark  
At a point $z$, open local connectivity is a strictly stronger
condition than simply local connectivity. However, the entire set $K$
is locally connected if and only if it is openly locally connected:
see Milnor~\cite[Section~16]{MiIntro}. We will discuss important
properties of locally connected sets in $\C$ in the appendix.

The following proposition will be the motor for many proofs of local
connectivity. 

\begin{proposition}[Trivial Fibers Yield Local Connectivity]
\label{PropFiberLocConn} \lineclear  
If a point of a compact connected full set $K\subset\C$ has a trivial
fiber, then $K$ is openly locally connected at this point. Moreover,
if the external ray at angle $\theta$ lands at a point $z$ with trivial
fiber, then for any sequence of external angles converging to $\theta$,
the corresponding impressions converge to $\{z\}$. In particular, if all
the fibers of $K$ are trivial, then $K$ is locally connected, all
external rays land, and the landing points depend continuously on the
angle.
\end{proposition}
\proof 
Consider a point $z\in K$ with trivial fiber. If $z$ is in the interior
of $K$, then $K$ is trivially openly locally connected at $z$.
Otherwise, let $U$ be an open neighborhood of $z$. By
Lemma~\ref{LemFibers}, any point $z'$ in $K-U$ can be separated from $z$
such that the separation avoids $z$ and $z'$. The region cut off from
$z$ is open; what is left is a neighborhood of $z$ having connected
intersection with $K$. By compactness of $K-U$, a finite number of such
cuts suffices to remove every point outside $U$, leaving another
neighborhood of $z$ intersecting $K$ in a connected set. Removing the
finitely many cut boundaries, an open neighborhood remains, and $K$ is
openly locally connected at $z$. Similarly, if $z$ is the landing point
of the $\theta$-ray, then external rays with angles sufficiently close
to $\theta$ will have their entire impressions in $\ovl U$ (although the
rays need not land). 
\qed

\remark 
The last statements of the proposition are always equivalent by
Cara\-th\'eo\-do\-ry's Theorem~\ref{ThmCaratheodory}. This is another
illustration of how closely fibers are related to Cara\-th\'eo\-do\-ry
theory. 

The converse to Proposition~\ref{PropFiberLocConn} is not true: local
connectivity at a certain point does not imply that the fiber of this
point is trivial. A counterexample is given in Figure~\ref{FigNonLocConn}.
However, if the set $Q$ in the definition of fibers is sufficiently big,
then local connectivity and triviality of fibers are equivalent for the
entire set. We give a general proof here, to be used for Julia sets in
Section~\ref{SecJulia}; for the Multibrot sets, there will be a direct
proof in \cite{FiberMandel}. (The following more local version of
this result seems plausible: whenever a point $z\in \partial K$ has a
neighborhood in $K$ such that $K$ is locally connected in this entire
neighborhood, then the fiber of $z$ is trivial for an appropriate
choice of $Q$.)

\begin{figure}[hbp]
%\vspace{60mm}
%\hspace{15mm}\special{picture NonLocConn}
\centerline{
\psfig{figure=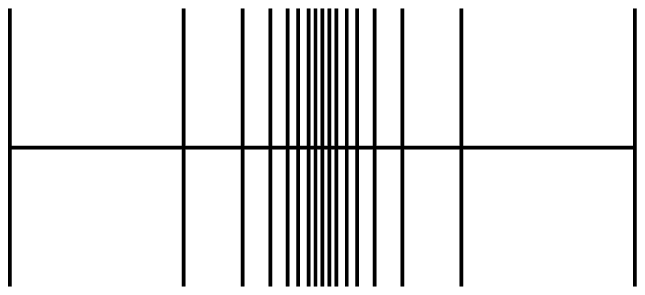,width=60mm}
}

\LabelCaption{FigNonLocConn}{A compact connected full set which is not
locally connected. It is locally connected at the center, but the fiber
of the center contains a vertical line segment, no matter which rays
are used to construct the fibers.} 
\end{figure}

\begin{proposition}[Local Connectivity Makes Fibers Trivial]
\label{PropLocConnFiber} \lineclear
Let $K\subset\C$ be a compact connected full set which is locally
connected. Suppose that $Q$ is a dense subset of\/ $\R/\Z$. Then all
fibers of $K$ are trivial provided that the following three conditions
are satisfied:
\begin{enumerate}
\item
whenever three external rays land at a common point, all their angles
are in $Q$; 
\item
if there exists an open interval $I$ of angles such that all the rays
with angles in $I$ land at different points, and each of their
landing points is also the landing point of some other ray, then there
exist angles $\theta\in I\cap Q$ and $\theta'\in Q$ such that the
corresponding rays land together. 
\item
if a point on the boundary of an interior component of\/ $K$ disconnects
$K$, then all the external rays of this point are in $Q$.
\end{enumerate}
\end{proposition}
\remark
In our applications, the second condition is usually void because the
landing points of the rays in $I$ would define an embedded arc in $K$
which contains no branches or decorations on at least one side. The
only quadratic polynomial where this condition applies is $z^2-2$ for
which the Julia set is an interval (for $z^2$, where the Julia set is
a circle, no two rays land together). In addition, the requirement in
this condition is usually satisfied anyway: in general, there is a
dense subset of external angles contained in $Q$ such that the
corresponding rays land together with another ray in $Q$. 

In the third condition, the disconnecting boundary point is the
landing point of at least two external rays by
Lemma~\ref{LemDisconnect}.
\proof
Since $K$ is locally connected, it is pathwise connected
by Lemma~\ref{LemLocConnPath}. Consider a connected component of the
interior of $K$ and let $Y$ be its closure. Then there is a dense subset
of $\partial Y$ (with respect to the topology of $\C$) consisting of
points which are landing points of rays in $Q$: if $U$ is an open set
intersecting $\partial Y$, then it either contains a boundary point of
$Y$ disconnecting $K$ (which is the landing point of a ray in $Q$ by the
third hypothesis), or density of $Q$ supplies a landing point of a ray
in $Q$ within $\partial Y\cap U$. Since local connectivity of $K$ is
equivalent to local connectivity of $\partial K$ (Carath\'eodory's
Theorem~\ref{ThmCaratheodory}), every boundary point of $Y$ is accessible
from the inside of $Y$. By Lemma~\ref{LemFibers}, the fiber of every
point in the interior of $Y$ is trivial, and every boundary point of $Y$
can be separated from any other point in $Y$ as well. Hence no fiber of
any point contains more than a single point on the closure of any
connected component of the interior.

Suppose that there is a fiber which is not trivial and denote it $Y$. It
has no interior, so we have $Y\subset\partial K$ by
Lemma~\ref{LemFibers}. Let $z_1\neq z_2$ be two points in $Y$ and let
$\gamma\subset\partial K$ be a simple closed curve connecting them; such
a curve exists by local connectivity of $\partial K$. We have
$\gamma\subset Y$ because otherwise $\gamma\cup Y$ would enclose an open
subset of $\C$ and thus an interior component of $K$, and the fiber $Y$
would meet more than a single point on the closure of this interior
component. Any point $z$ on the interior of $\gamma$ is the landing
point of at least two external rays, one from either side of $\gamma$:
this is because the curve $\gamma$ cuts every sufficiently small disk
$D$ around $z$ in at least two parts, and both parts must intersect the
exterior of $K$. No interior point of $\gamma$ can be the landing point
of three or more external rays because otherwise we could  separate
$\gamma$ and thus $Y$. Let $\alpha$ be an external angle of
$z$; then rays at angles near $\alpha$ must land near $z$, and if they
did not land on $\gamma$, then $\gamma$ would have a branch point near
$z$. Therefore, rays with angles sufficiently close to $\alpha$ land at
interior points of $\gamma$, and by the second assumption, some of them
must be in $Q$ and landing together with another ray in $Q$. This ray
must come in from the other side of $\gamma$, and we can separate $Y$
again.

It follows that every fiber of a locally connected set $K$ is trivial,
provided that $Q$ is sufficiently large so as to satisfy the stated
conditions. 
\qed

\remark
For any compact connected and full $K\subset\C$ which is locally
connected, there is always a countable collection $Q$ of external angles
for which all the fiber become trivial: The second condition requires
only countably many rays. The first and third conditions specify
countably many points: the number of branch points is countable by
Lemma~\ref{LemBranchCountable} below; similarly, the number of interior
components is obviously countable, and each has at most countably many
disconnecting boundary points by Corollary~\ref{CorProjectCount}.
The problem is that some of these points might be the landing points of
uncountably many external rays. Even in that case, the number of
connected components any such point disconnects $K$ into is countable by
Lemma~\ref{LemBranchCount}, and countably many rays at every branch
point suffice to separate any two of the connected components of the
complement. The proposition remains true with these weakened hypotheses.
However, in our applications only finitely many rays land at any single
point and the given form of the proposition suffices.

\begin{lemma}[Branch Points Countable]
\label{LemBranchCountable} \lineclear
For any compact connected and full subset $K$ of\/ $\C$ and any
$\eps>0$, the number of points which are the landing points of at
least three external rays with mutual distance at least $\eps$ is
finite and bounded above independently of $K$. In particular, the
number of points which are the landing points of at least three rays
is countable. 
\end{lemma}
\remark
The ``distance between external angles'' will be the distance between
their external angles in $\R/\Z$, so that the maximal distance is
$1/2$.
\proof
We will follow a suggestion of Saeed Zakeri. Parametrize the boundary of
$\disk$ by external angles in $\R/\Z$. When three external rays land
at a common point, mark this by a Euclidean triangle in $\disk$ with
vertices at the boundary points of $\disk$ corresponding to the external
angles of the rays. Triples of rays landing at distinct points will then
give rise to disjoint triangles. If all the angles of the triangle have
mutual distance at least $\eps>0$, then the Euclidean area of the
triangle will be bounded below. Since the total area of the disk is
finite, the number of such triangles is finite. The second claim follows.
\qed

\newsection{Fibers of Filled-in Julia Sets}
\label{SecJulia}

In this section, we will apply the general concept of fibers from
Section~\ref{SecFibers} to connected filled-in Julia sets, where the set
$Q$ of external angles will always contain the rational numbers
$\Q/\Z$ and sometimes countably many further angles. We will always
assume the set $Q$ to be forward and backward invariant under
multiplication by the degree, so that the set of corresponding dynamic
rays is invariant under the dynamics. Several of the results in this
section will be valid for arbitrary connected Julia sets of polynomials
(which we will then always assume to be monic), while others are proved
only for unicritical polynomials. 

First we show that every bounded Fatou component has zero or
infinitely many boundary points which are accessible from inside and
outside, which makes the relation ``is in the fiber of'' symmetric. We
will then discuss branch points of unicritical Julia sets: the analog
to the Branch Theorem for the Multibrot sets
\cite[Theorem~\MandelBranch]{FiberMandel} is Thurston's No Wandering
Triangles Theorem~\ref{ThmThurstonTriangles}. We know that local
connectivity and triviality of all fibers are equivalent for some
choice $Q$ of external angles. The No Wandering Triangles Theorem will
allow to specify the set $Q$. 

For the Multibrot sets, local connectivity implies that every
connected component of the interior is hyperbolic; similarly, a
corollary to Thurston's theorem is that locally connected Julia sets
of unicritical polynomials do not have wandering domains, i.e., all
their Fatou components  are eventually periodic. This result holds for
arbitrary rational maps by Sullivan's Theorem, and we will assume it
throughout. 

Finally, we will establish the ``nice'' situation of
Lemma~\ref{LemFibersNice} for certain Julia sets: landing points of
rational rays have trivial fibers, so the filled-in Julia sets split
into equivalence classes of points with coinciding fibers. The
analogous statement for the Multibrot sets is discussed in
\cite{FiberMandel}.

\begin{lemma}[Accessibility of Interior Components of Julia Sets]
\label{LemFibersJulia} \lineclear
Consider an arbitrary polynomial with connected filled-in Julia set. 
Every bounded Fatou component corresponding to an attracting or
rationally indifferent periodic orbit has infinitely many boundary
points which are landing points of dynamic rays at rational angles
and which are also accessible from within the component. 

Whenever any bounded Fatou component eventually lands on a periodic
orbit of Siegel disks and has a single boundary point which is
accessible from inside and which is also the landing point of an
external dynamic ray, then every Fatou component on the same grand
orbit has countably many such boundary points. This always happens when
the Julia set is locally connected. However, the corresponding external
angles are in no case rational. 
\end{lemma}
\proof
Denote the filled-in Julia set by $K$ and consider a bounded periodic
Fatou component. If this Fatou component belongs to an attracting or
rationally indifferent orbit, there is at least one boundary point which
is fixed under the first return map of the Fatou component. 
This point must be repelling or rationally indifferent. It is thus
the landing point of at least one rational dynamic ray, and it is
accessible from within its Fatou component. Every Fatou component which
eventually maps onto $U$ then inherits countably many points on its
boundary which are all accessible from inside and which are landing
points of rational dynamic rays.

The only further type of Fatou components of polynomials are Siegel
disks and their preimages. If a boundary point $z$ of a Siegel disk is
accessible both from inside and outside, it is the landing point of an
external and of an internal ray (by Lindel\"of's
Theorem~\ref{ThmLindeloef}; we define internal rays with respect to the
periodic point at the center as the base point). Since the dynamics on
internal rays is an irrational rotation, the point $z$ cannot be
periodic: otherwise, it would be the landing point of two (and even
countably many) internal rays, and the region between them would have to
be contained entirely within the Siegel disk because the filled-in Julia
set is full and the boundary of the Siegel disk is contained in the
boundary of the filled-in Julia set. But then an open interval of
internal angles would have to land at the same point, which is a
contradiction to the theorem of the Riesz brothers
\cite[Theorem~A.3]{MiIntro}. Any boundary point accessible from inside
and outside thus gives rise to countably many such points, and their
external angles are all irrational. Again, every Fatou component which
eventually maps onto this Siegel disk inherits countably many boundary
points with the specified property.  

If the filled-in Julia set is locally connected, then there are many
such boundary points: any boundary point of the Siegel disk which is
accessible from inside will do the job, and these are dense (in fact,
by Lemma~\ref{LemLocConnInterior}, the boundary of the Siegel disk itself
is locally connected, and each of its boundary points is accessible from
inside and outside).
\qed

The following corollary shows that the relation ``is in the fiber of'' is
symmetric for arbitrary connected Julia sets. 

\begin{corollary}[Fibers are Symmetric]
\label{CorFiberSep} \lineclear 
Consider an arbitrary polynomial with connected filled-in Julia set\/
$K$. Define fibers of\/ $K$ using an arbitrary choice of the set\/ $Q$
of external angles which is forward and backward invariant (subject to
the usual two conditions that $Q$ be countable and that all rays with
angles in $Q$ actually land). Let $z,z'\in K$ be two points such that
$z'$ is not in the fiber of\/ $z$. Then there is a separation line
separating $z$ and $z'$ which avoids these two points, and $z$ is not in
the fiber of $z'$.
\end{corollary}
\proof
By Lemma~\ref{LemFibers}, the claim can fail only if $z$ is in the
interior of $K$ and the connected component of the interior of $K$ which
contains $z$ has exactly two boundary points which are accessible from
inside and outside. But every Fatou component will eventually map onto a
periodic Fatou component corresponding to an attracting or rationally
indifferent periodic point or onto a Siegel disk. For those Fatou
components, the number of boundary points which are accessible from
inside and outside is either zero or infinite. 
\qed

The principal goal in this section is to specify a choice $Q$ of
external angles for which the fibers of a locally connected
unicritical Julia set are trivial. We have to check three conditions
in Proposition~\ref{PropLocConnFiber}: the first one is easy to
satisfy and  the second one is usually void. For the third condition,
we need a theorem due to Thurston~\cite[Theorem~II.5.2]{Th} which is
still unpublished. It is the dynamic analog to the Branch
Theorem~\cite[Theorem~\MandelBranch]{FiberMandel} for the Multibrot
sets, stating that branch points have rational external angles.
Thurston states his theorem only for quadratic polynomials, but his
proof works for all unicritical polynomials. With his permission, we
give his proof here. It is slightly modified using an idea of Saeed
Zakeri.

\begin{theorem}[No Wandering Triangles]
\label{ThmThurstonTriangles} \lineclear
If three dynamic rays of a unicritical polynomial with connected
Julia set land at a common point, then the landing point is either
periodic or preperiodic, or it eventually maps through a critical
point.
\end{theorem}
\remark
If the landing point is on a repelling or rationally indifferent
orbit, then the rays are all periodic or all preperiodic and have
thus rational angles. The only other conceivable case is that the
landing point is a Cremer point and all the rays landing there are
irrational. As far as I know, it is not known whether that can
possibly happen.

Thurston proves his theorem in an abstract setting using
``laminations'', related to the pinched disk model of the Julia set.
That way, he does not have to worry whether certain dynamic rays land
at all. We will use the theorem only for Julia sets which are locally
connected, so all dynamic rays land and there is no Cremer point. 

\proof 
External angles are parametrized by $\Circle=\R/\Z$; identify this
set with $\partial\disk$. Assume that three dynamic rays at
angles $\theta_1,\theta_2,\theta_3$ land at a common point. If the
theorem is false, then the forward orbit never repeats and never maps
through the critical point. We will suppose that in the following. The
three angles are necessarily irrational and will remain distinct under
forward iteration. For every $k\geq 0$, the dynamic rays at angles
$d^k\theta_1, d^k\theta_2, d^k\theta_3$ also land at a common point.

On $\partial\disk$, connect the three points $\theta_i$ pairwise
by Euclidean straight lines, yielding a Euclidean triangle in $\disk$
which represents the landing point of these three rays: every side of
the triangle stands for a ray pair. Since ray pairs landing at
different points do not cross and all the landing points are
different, we obtain an infinite sequence of disjoint image triangles
connecting the angles $d^k\theta_1,d^k\theta_2,d^k\theta_3$: a
wandering triangle.

Because of the $d$-fold rotation symmetry of the Julia sets, every
triangle has $d-1$ rotated counterparts, and adding these in still
leaves the triangles non-intersecting: each of these extra triangles
corresponds to the landing point of three rays which maps in one step
onto the orbit of the initial triangle. This is where we are using the
assumption that the polynomials are unicritical. 

We will measure the lengths of a triangle side (i.e., of a ray pair)
as usual as the unsigned distance along $\Circle$ between the
corresponding angles. The maximum distance between any two points is
therefore $1/2$, realized for points straight across. In fact,
because of the rotation symmetry and since triangle sides never
cross, no side can have length $1/d$ or more (except for rays landing
at the critical point, which is the center of symmetry; in that case,
we discard the initial triangle and consider only the remaining
orbit). If a side of a triangle has length $s<1/d$, then after
multiplication by $d$, the image side will have length
$\min\{ds,1-ds\}$ (measuring the short way around the circle), so
that sides with lengths less than $1/(d+1)$ will be mapped to longer
sides, while those with lengths greater than $1/(d+1)$ will shrink in
length. Short sides of length $\eps$ are images of sides of length
$\eps/d$ or of length $1/d-\eps/d$, so they are images of very short
or of very long sides.

By Lemma~\ref{LemBranchCountable}, there can be only finitely many
points which are landing points of three dynamic rays with mutual
distance at least $\eps$, for any $\eps>0$. Therefore, if there is a
wandering triangle, then the lengths of the respective shortest sides
must converge to zero. It follows that there can be no upper bound
less than $1/d$ for the lengths of sides because a new shortest side
can be the image only of a very long side. Therefore, there exists a
sequence $k_1,k_2,\ldots$ of iteration steps such that the longest
side of the $k_1$-th image of the wandering triangle has length
$l_1>1/(d+1)$ and the image after $k_{i+1}$ steps has a longest side
of length $l_{i+1}>l_i$. Denote the respective triangles by $T_i$ and
denote the lengths of its other two sides by $l'_i$ and $l''_i$ such
that $l_i\geq l'_i\geq l''_i$. We want this sequence to be maximal in
the following sense: the first image of $T_i$ with a side of length
exceeding $l_i$ is already $T_{i+1}$.

The side of $T_i$ with length $l_i$ and its $d-1$ symmetric rotates
cut the disk into $d+1$ pieces, of which one contains the origin and
is rotation symmetric. Denote this piece by $C_i$. Since
$l_{i+1}>l_i>1/(d+1)$, the side with length $l_{i+1}$ must be
contained in $C_i$, together with the triangle $T_{i+1}$ it belongs
to. Therefore, we also have $l'_{i+1}>l_i$: two sides of a new
triangle will be longer than the longest side of an old triangle. It
follows that $l_{i+1}>l'_{i+1}>l''_{i+1}$ with strict inequality; this
holds for every $i$.

We claim that the two long sides of any triangle $T_i$ will, after
$k_{i+1}-k_i$ iterations, map onto the two long sides of $T_{i+1}$.
Indeed, the shortest side of $T_{i+1}$ has length less than $1/d-l_i$
because this is the length of the intervals in which $C_i$ meets
$\Circle$. However, the image of the longest side of $T_i$ has length
$1-dl_i=d(1/d-l_i)$, so it is already too long for the shortest side
of $T_{i+1}$; the image of the middle side of $T_i$ is even longer.
If the two long sides of $T_i$ want to become shorter, they must
first be longer. The first time that this happens they are on the
triangle $T_{i+1}$, proving the claim.

Perhaps not unexpectedly, we obtain a contradiction by looking at the
orbit of the shortest sides, which must always map to the shortest
sides. No matter how short it started, it will eventually have length
at least $1/(d+1)$ and might then get shorter. But in order to map to
the shortest side of a triangle $T_{i+1}$, it must have been very
short in $T_i$ or longer than $l_i$. The second option is clearly
impossible, and the first can happen only a finite number of times.
To acquire a new shortest length, it must have been very long before,
and that happens only at the $T_i$. Here is the contradiction.
\qed

\remark
This theorem has recently been generalized by Kiwi~\cite{Kiwi} to
arbitrary polynomials with connected Julia sets: he has a ``No Wandering
Polygon'' Theorem, but the number of sides of his polygons depends on
the degree.

\begin{lemma}[Dynamics of Fibers]
\label{LemFiberDynamics} \lineclear
For any polynomial with connected Julia set and any choice of the set\/
$Q$ which is forward and backward invariant, the dynamics maps the
fiber of any point as a possibly branched cover onto the fiber of the
image point.
\end{lemma}
\proof
Let $K$ be the filled-in Julia set of a polynomial $p$ having degree
$d$. We know from Corollary~\ref{CorFiberSep} that whenever one point
is not in the fiber of another, then these two points can be separated
by a separation line avoiding both points. 

Choose a point $z\in K$, let $z':=p(z)$, let $Y'$ be the fiber of $z'$
and let $Y_0$ be the connected component of $p^{-1}(Y')$ containing $z$.
Since $p\colon\C\to\C$ is a branched covering, its restriction to
$Y_0$ is also a branched covering. Denote the fiber of $z$ by $Y$. We
will first show that $Y\subset Y_0$. 

Choose an arbitrary point $z'_1\in K-Y'$. Then there is a separation
line separating $z'$ and $z$. If this separation line does not contain a
critical value, then its pull-back under $p$ will be $d$ separation
lines, and every inverse image of $z'_1$ can be separated from $z$ by
one of them. If the separation line does contain a critical value, then
it is still possible to separate every inverse image of $z'_1$ from $z$
by a separation line made up of parts of the inverse images of the given
separation line. Therefore, $Y\subset Y_0$. 

If already $Y=Y_0$, then we are done. If not, let $z_1$ be a point in
$Y_0-Y$ and consider a separation line $\gamma$ between $z$ and $z_1$.
If it is a ray pair which does not land at a critical point of $p$,
then the image of $\gamma$ is another separation line. Since $\gamma$
runs through $Y_0$ and separates it, its image will run through $Y$ and
separate it. This is impossible. If $\gamma$ is a separation line
running through an interior component $U$ of $K$ and the images of its
two dynamic rays are different, then the image of $\gamma$ is again a
separation line, possibly after modifying it within $p(U)$ so that the
new separation traverses $p(U)$ in a simple curve. All the fibers of
points in $U$ and $p(U)$ will then be trivial. Since $\gamma$
disconnects $Y_0$, which is connected, the landing point of at least one
of the two dynamic rays in $\gamma$ will have a neighborhood in $Y_0$
which is disconnected by $\gamma$. The new separation line will then
separate $Y'$ at the image point, which is again impossible. 

Therefore, if $Y\neq Y_0$, then any separation line $\gamma$ which
separates $Y_0$ has the property that its two dynamic rays have the
same image rays, or it is a ray pair landing at a critical point. If
$\gamma$ runs through an interior component $U$ of $K$, then there are
countably many further dynamic rays landing at $U$ which are accessible
from inside, and it is easy to manufacture a new separation line which
still separates $Y_0$ but which will not collapse when mapped forward,
so the argument above applies: an impossibility again. The last case is
that $\gamma$ is a ray pair landing at a critical point. Removing from
$Y_0$ the part which is separated from $z$, it is easy to check that
$p$ induces a covering from the rest onto $Y$: the only place where we
have to check this is at the landing point of the ray pair, and there is
no problem. Since there are only finitely many critical points, and
these have only finitely many rays landing, there are only finitely
many such ray pairs. After finitely many cuts in $Y_0$, we obtain the
fiber $Y$, and $p\colon Y\to Y'$ is a branched covering. If there are
branch points at all, these are critical points of $p$.
\qed

We will now show that, at least for many Julia sets, the landing
points of rational rays have trivial fibers. The corresponding
statement for Multibrot sets can be found in \cite{FiberMandel}.

\begin{theorem}[Repelling Periodic Points Have Trivial Fibers]
\label{ThmPeriodicFiber} \lineclear 
Consider a polynomial with connected filled-in Julia set and define its
fibers for $Q=\Q/\Z$, together with the grand orbits of all the rays
landing at those critical values which are on the boundary of periodic
Siegel disks (if any). Let $z$ be a repelling periodic or preperiodic
point and suppose that all the points on its forward orbit can be
separated from all the critical values and from all the points on
closures of periodic bounded Fatou components. Then the fiber of $z$ is
trivial.
\end{theorem}
\proof
By Lemma~\ref{LemFiberDynamics}, the fiber of any point is trivial
whenever it every maps to a point with trivial fiber. Therefore, we may
assume that $z$ is periodic. By switching to an iterate, we may assume
$z$ to be a fixed point. Denote the corresponding (iterated) polynomial
by $p$ and let $K$ be its filled-in Julia set.

Every point on the closure of a periodic bounded Fatou component can be
separated from $z$ by a separation line. Since every such line separates
from $z$ an open subset of the closure of this periodic Fatou component,
a finite number of separation lines suffices to separate the entire
closure of this Fatou component (in fact, a single line will do the
job). The total number of periodic Fatou components is finite, so there
is a finite number of separation lines separating $z$ from all the
bounded periodic Fatou components and from all the critical values.
Denote this collection of ray pairs by $S_0$ and let $U_0$ be the
neighborhood of $z$ which is not separated from $z$ by separation lines
in $S_0$.

Consider all the separation lines in $S_0$ which are not ray pairs. They
will then traverse bounded Fatou components, so all but finitely many
of their images under forward iteration will intersect bounded periodic
Fatou components. Therefore, only finitely many of these forward images
can intersect and cut $U_0$, and none of them can meet $z$. A similar
argument applies to those ray pairs in $S_0$ which have irrational
external angles, so they necessarily land on the boundary of periodic or
preperiodic Siegel disks. Let $U_1$ be the connected component of $z$ in
$U_0$ minus these finitely many separation lines. 

Now we look at separation lines bounding $U_1$ which are ray pairs at
rational angles. Their landing points are periodic or preperiodic. Then
all these separating ray pairs have finite forward orbits. Consider all
the finitely many ray pairs on these forward orbits, except those
landing at $z$. They might possibly disconnect $U_1$. Let $U_2$ be
the connected component of $z$ in $U_1$ minus these finitely many ray
pairs. Consider an arbitrary equipotential of $K$ and let $U$ be the
subset of $U_2$ within this equipotential. Then $p$ restricted to $U$ is
a conformal isomorphism onto its image, and $p$ cannot send boundary
points of $U$ into the interior of $U$. 

Since $U$ is full and contains no critical point, the branch of $p^{-1}$
fixing $z$ can be extended throughout $U$. All the ray pairs and
separation lines bounding $U$ are mapped into $U$ or to its boundary:
if they are mapped outside of $U$, then a separation line on the
boundary of $U$ is inside $p^{-1}(U)$, and mapping $p^{-1}(U)$ forward
under $p$ sends a bounding ray pair into $U$, which we had excluded
above. Since the equipotential bounding $U$ is mapped to a lower
equipotential under $p^{-1}$, the branch of $p^{-1}$ fixing $z$ maps $U$
into itself. 

Therefore, the restriction of $p^{-1}$ to $U$ is a holomorphic self-map
of $U$ with an attracting fixed point at $z$. Each of the finitely many
separation lines bounding $U$ is either mapped eventually into $U$, or
it is periodic. The latter case is impossible because the separation
line would necessarily have to be a ray pair at rational angles, all
parabolic periodic points are separated from $z$ by assumption, and
repelling periodic points would have to attract nearby points under
iteration of $p^{-1}(U)$, while the interior of $U$ has to converge to
$z$ by Schwarz' Lemma. 

Therefore, all of $\ovl U$ converges to $z$ under iteration of
$p^{-1}(U)$. For every $\eps>0$, there is an $n$ such that
$p^{\circ(-n)}(\ovl U)$ is contained in the $\eps$-neighborhood of $z$.
But that means that no point $z'\in K$ with $|z'-z|>\eps$ can be in the
fiber of $z$. Since $\eps$ was arbitrary, the fiber of $z$ is trivial.
\qed

\remark
It is important to require that $z$ can be separated from closures of
periodic Siegel disks. The separation from other periodic bounded Fatou
components (attracting or parabolic) is for convenience and does not
seem essential. Similarly, a related proof will probably transfer the
proof from repelling to parabolic periodic points. For unicritical
polynomials, the presence of attracting or parabolic orbits makes all
the fibers of the Julia set trivial anyway.

From now on, we will restrict to filled-in Julia sets of unicritical
polynomials. For these, we can now specify a set $Q$ for which
triviality of all fibers is equivalent to local connectivity of the
Julia set. We already know from Proposition~\ref{PropFiberLocConn} that
triviality of fibers implies local connectivity, so we only state the
converse. 

\begin{proposition}[Locally Connected Julia Sets have Trivial Fibers]
\label{PropLocConnJuliaFiber} \lineclear
If the filled-in Julia set of a unicritical polynomial is locally
connected, then all its fibers are trivial for the choice $Q=\Q/\Z$
unless there is a Siegel disk; in that case, all fibers are trivial when
$Q=\Q/\Z$ together with the grand orbits of the angles of all the rays
landing at the critical value.
\end{proposition}
\remark
A locally connected Julia set of a polynomial can never have a Cremer
point; see Milnor~\cite[Corollary~\MiCor]{MiIntro}. In the case of a
Siegel disk, all the rays we really need are the rays in $\Q/\Z$ and
those landing at the critical point and on its backwards orbit; the
extra rays are just taken in to have invariance of the rays in $Q$ under
the dynamics. We will see below that a single ray lands at the critical
value and at every point of its forward orbit. The separation lines
through periodic Siegel disks which we can obtain from such rays can be
replaced by lines through precritical points.
\proof
The No Wandering Triangles Theorem implies that three or more rays
landing at a common point either have rational angles, or the landing
point eventually maps through the critical point. 

First we discuss the case that the filled-in Julia set has no
interior. Being locally connected, it is a {\em dendrite:} any pair of
points can be connected by a unique arc within the Julia set
(Lemma~\ref{LemLocConnPath}). Separation lines are just ray pairs at
rational angles.

The critical point cuts the Julia set into two parts, to be labelled
$\0$ and $\1$, and this partition defines a symbolic itinerary for any
point which is not a pre-critical point. The subset of the Julia set
with identical first $k$ entries in the itinerary is connected, and no
two points have identical itineraries forever (otherwise, an entire
interval of external angles would have to have the same itinerary).
Therefore, precritical points are dense on any subarc of the Julia
set. 

Within the dendrite Julia set, the critical orbit spans an invariant
subtree (a postcritically infinite Hubbard tree), and the critical
value is an endpoint of this tree. It follows that the critical point
cannot be a branch point of the Hubbard tree, so all its branch points
are periodic or preperiodic. 

The critical value is a limit point of periodic points in the tree: if
$z_n$ is a precritical point on the tree such that the interval
between $z_n$ and $c$ contains no point which maps before $z_n$ onto
the critical point, then there is a homeomorphic forward image of the
interval $[z_n,c]$ which maps $z_n$ onto $c$, producing a periodic
point on this interval. By density of precritical points, periodic
points are dense on every subarc of the Julia set. It follows that
any two given points in the Julia set can be separated by a periodic
point, which is necessarily repelling, and the rays landing at this
periodic point separate the two given points. Therefore, all fibers
are trivial, even when $Q$ only contains periodic angles.

We now consider the case that the filled-in Julia set has interior. 
We will prove the result by checking the conditions in
Proposition~\ref{PropLocConnFiber}. 

If the bounded Fatou components correspond to an attracting or
parabolic orbit, then the critical orbit is in the Fatou set and
$Q=\Q/\Z$; otherwise, we have a Siegel disk and the critical point is
in the Julia set. In that case, $Q$ contains countably many further
rays. In both cases, all the external angles of branch points are in
$Q$ by the No Wandering Triangles Theorem, and the first condition of
Proposition~\ref{PropLocConnFiber} is always satisfied. Moreover, the
number of rays landing at any given point is well known to be finite.

If there is an open interval of external angles of length less than
$1/d$ such that all the corresponding dynamic rays land at different
points, then multiplication by $d$ yields another longer interval with
the same property. Restricting to a subinterval of length $1/d^n$ for
an appropriate integer $n$ and iterating this argument, it follows
that all dynamic rays land at different points. The second condition
is thus always void.

We are assuming that there is a periodic cycle of bounded Fatou
components. Let $U$ be one such component and let $z_1$ be a boundary
point of $U$ which disconnects the filled-in Julia set. Then at
least two dynamic rays land at $z_1$ by Lemma~\ref{LemDisconnect},
but the total number of rays at $z_1$ is always finite. Let
$\theta_1$ and $\theta'_1$ be the angles of two rays landing at $z_1$
so that they separate as much as possible from $U$. Denote the period
of $U$ by $n$. Iterating the $n$-th iterate of the polynomial, we
obtain a sequence $z_2,z_3,\ldots$ of boundary points of $U$ and two
sequences of dynamic rays at angles $\theta_2,\theta_3,\ldots$ and
$\theta'_2,\theta'_3,\ldots$. 

Each ray pair $(\theta_k,\theta'_k)$ cuts away an open interval of
external angles from $U$, so that the projection (as defined after
Lemma~\ref{LemProjectComp}) of external rays within such an interval
yields the point $z_k$. If all the points $z_k$ are different and all
$\theta_k\neq\theta'_k$, then they will cut away infinitely many
intervals which must all be disjoint. Therefore, their lengths must
shrink to zero. However, when such intervals are short, then their
lengths are multiplied by the degree $d$ of the polynomial in every
step and by $d^n$ under the first return map of $U$, so there will
always be intervals with lengths bounded below. This is a
contradiction.

Therefore, either the point $z_1$ is periodic or preperiodic and then
its external angles are rational, or it is in the backwards orbit of
the critical point. In both cases, its external angles are in $Q$,
satisfying the third condition of Proposition~\ref{PropLocConnFiber}
and finishing the proof also in the case when there are bounded Fatou
components. 
\qed

\remark
The second proof given also applies when there are no bounded Fatou
components, but it gives a weaker result because it specifies a larger
choice of $Q$.

\hide{
If there are no bounded Fatou components, the proof above still shows
that local connectivity implies that fibers are trivial when $Q$ is the
set of rational rays, together with the grand orbits of all the rays
landing at the critical value. However, we want to prove that the Julia
set has trivial fibers with $Q=\Q/\Z$. Before doing that,}

We now state some observations which came out of the proof. They are
all known.

\begin{corollary}[Disconnecting Boundary Points of Fatou Components]
\label{CorDisconnectJulia} \lineclear
Let $z$ be a boundary point of a bounded Fatou component of a
unicritical polynomial and assume that it disconnects the Julia set.
If the Fatou component corresponds to an attracting or parabolic
periodic orbit (in which case the Julia set is known to be locally
connected), then $z$ is a periodic or preperiodic point. If the Fatou
component corresponds to a Siegel disk, and we assume the Julia set
to be locally connected, then $z$ will eventually map to the critical
point. In particular, the critical point of a unicritical polynomial
with a locally connected Julia set featuring a Siegel disk is on the
boundary of one of the periodic components of the Siegel disk, and
the critical value is the landing point of a unique dynamic ray.
\qedd
\end{corollary}

\begin{lemma}[Critical Point in Periodic Fiber]
\label{LemCriticalFiberPeriodic}\lineclear
Consider a unicritical polynomial and set $Q=\R/\Z$. If the fiber
containing the critical point is periodic of some period $n$, then the
polynomial is $n$-renormalizable and the critical fiber contains an
indifferent or superattracting periodic point of period $n$.
\end{lemma}
\proof
The statement is void or trivial if the Julia set is locally connected,
so we can in particular exclude hyperbolic or parabolic Julia sets. For
other parameters on the closure of the main hyperbolic component of a
Multibrot set, no two rational dynamic rays land together, and the
entire Julia set is a single fiber. The claim holds trivially for $n=1$. 
Otherwise, there is a unique repelling fixed point which is the landing
point of at least two dynamic rays. Denote this fixed point by $\alpha$.
The rays landing at $\alpha$ separate the critical point from the
critical value. If the critical fiber is periodic, its period must be at
least two.

Let $Q'\subset Q$ be the union of the rays at $\alpha$ together with
their entire backwards orbits. These are the rays usually used in the
construction of the Yoccoz puzzle. The critical fiber corresponding to
these rays will still be periodic of some period $n'$ dividing $n$,
again with $n'\geq$. It is quite easy to see and well known that the
polynomial is now $n'$-renormalizable (see e.g.\ Milnor~\cite[Lemma~2]
{MiLocConn}). After $n'$-renormalization, we have a new unicritical
polynomial with equal degree, and the critical fiber is still periodic
of period $n/n'$. If we are now on the closure of the main hyperbolic
component of the Multibrot set, the entire Julia set is a single fiber,
the critical fiber has period $1$ and contains a non-repelling fixed
point, and if it is attracting, then all fibers are trivial. It follows
that $n=n'$. For the original polynomial, the critical fiber must
contain an indifferent or superattracting periodic point of period
dividing $n$. 

If the renormalized polynomial is not on the closure of the main
hyperbolic component, then the period of the critical fiber is again at
least $2$, and we can repeat the argument. Since the period of the
critical fiber is reduced in every step, we must land after finitely
many steps on the closure of the main hyperbolic component.
\qed

\hide{
\remark
If a unicritical polynomial has a repelling periodic point which contains
the critical point in its fiber (for $Q=\Q/\Z$), then the fiber of the
critical point is periodic with some period $n$ and the polynomial turns
out to be $n$-renormalizable. \hide{\cite{FiberTuning}}   After
renormalization, we  have a repelling fixed point which contains the
critical point in its fiber. But it is easy to see that for parameters
outside of the main component of a Multibrot set or at rational boundary
points of the main component, all the fixed points can be separated from
the critical value and the critical point. Therefore, the only case where
the previous lemma does not apply is when there is an irrationally
indifferent periodic point. 
We see that, for all unicritical polynomials which have connected
filled-in Julia sets without indifferent periodic points, the fibers with
respect to $Q=\Q/\Z$ are ``nice'' (in the sense of
Lemma~\ref{LemFibersNice}), so that two points have identical fibers if
and only if their fibers intersect; fibers are thus an equivalence
relation on the filled-in Julia set. This applies in particular to all
infinitely renormalizable polynomials. 
}

\begin{corollary}[Impressions of Rational Dynamic Rays]
\label{CorRatRayImpress} \lineclear
For any unicritical polynomial in $\M_d$ without indifferent periodic
points, the impression of any dynamic ray at a rational angle is
always a single point. Fibers of any two points (for $Q=\Q/\Z$) are
either disjoint or equal and have the ``nice'' property of
Lemma~\ref{LemFibersNice}.
\end{corollary}
\proof 
If the critical fiber is periodic, then there is either a
superattracting or an indifferent orbit by
Lemma~\ref{LemCriticalFiberPeriodic}. The indifferent case is
excluded. If there is a superattracting orbit, or if the critical
fiber is not periodic, then every repelling periodic point can be
separated from the critical value. By Theorem~\ref{ThmPeriodicFiber},
the fibers of repelling periodic points are trivial, and they contain
the entire impressions of all the rays landing there by
Lemma~\ref{LemImpressFiber}. 

This is obvious since rational dynamic rays always land, and the
impression of any ray is contained in the fiber of its landing point
by Lemma~\ref{LemImpressFiber}. This establishes the ``nice'' situation
of Lemma~\ref{LemFibersNice}, and the Julia set splits into equivalence
classes of points having intersecting and thus identical fibers.
\qed
\remark
For non-infinitely renormalizable quadratic polynomials, this is a
special case of a theorem of Yoccoz~\cite[Theorem~II]{HY}. Very
recently, J.~Kiwi~\cite{Kiwi} has independently proved this theorem for
arbitrary polynomials with connected Julia sets and with all periodic
points repelling.

\appendix
\newsection{Compact Connected Full Sets in the Plane}
\label{SecLocConn}

In this appendix, we will discuss compact connected full (and
sometimes locally connected) subsets in $\C$ and describe certain
properties which we will need in the main text. Local connectivity has
been defined in Definition~\ref{DefLocConn}.

Of principal importance is that local connectivity implies pathwise
connectivity, i.e., any two points can be connected by a continuous
image of an interval. In fact, we can connect them by a homeomorphic
image of an interval, a property known as arcwise connectivity.
\begin{lemma}[Local Connectivity Implies Arcwise Connectivity]
\label{LemLocConnPath} \lineclear
Every compact connected and locally connected subset of\/ $\C$ is
arcwise connected and locally arcwise connected.
\qedd
\end{lemma}
For a proof, see Douady and Hubbard~\cite[Expos\'e~II]{Orsay},
or Milnor~\cite[Section~16]{MiIntro}.

Another important result is Carath\'eodory's Theorem, which is also
described in \cite{Orsay} and \cite{MiIntro}.
\pagebreak[2]
\begin{theorem}[Carath\'eodory's Theorem]
\label{ThmCaratheodory} \lineclear
Let $K$ be a compact connected and full subset of\/ $\C$. Then $K$
is locally connected if and only if\/ $\partial K$ is locally connected,
or if and only if all the external rays of\/ $K$ land with the landing
points depending continuously on the external angles. In that case,
every boundary point is the landing point of at least one external
angle.
\qedd
\end{theorem}
If $K$ is locally connected, then the map from external angles to
the corresponding landing points is known as the {\em Carath\'eodory
loop} of $K$, and it is surjective onto $\partial K$.

\begin{lemma}[Interior Component Locally Connected]
\label{LemLocConnInterior} \lineclear
Consider a compact, connected and full subset of\/ $\C$ which is
locally connected. Then any connected component of the interior has
locally connected boundary.
\end{lemma}
\proof
Denote the original locally connected set by $K$, let $U_0$ be an
interior component of $K$ and let $K_0$ be the closure of $U_0$. Let
$z\in \partial K_0$ and let $U$ be an open neighborhood of $K_0$. By
local connectivity of $K$, there is a neighborhood $V\subset U$ of $z$
such that $V\cap K$ is connected. We claim that $V\cap K_0$ is connected.

Suppose not. Then let $K_1$ and $K_2$ be two connected components of
$V\cap K_0$. Since $V$ is open, both $K_1$ and $K_2$ contain interior
points of $K_0$, and there is a curve $\gamma$ in the interior of $K_0$
connecting $K_1$ and $K_2$. Obviously, this curve cannot be contained
entirely within $V$. Since it connects two points in $V\cap K$, the set
$V\cap K$ and $\gamma$ together disconnect $\C-(V\cap K)$. Let $W$ be
open subset of $\C-(V\cap K)$ which is disconnected from $\infty$ by
the curve $\gamma$. Since $\partial W\subset K$ and $K$ is full, we have
$W\subset K$. And since the interior component $U_0$ intersects $W$, it
follows $W\subset U_0$. The entire boundary of $W$, except the part on
$\gamma$, is contained in $V$ by construction, and it is also contained
in $K_0$. Therefore, $K_1$ and $K_2$ are connected within $V\cap K_0$,
contrary to our assumption.
\qed

\remark
It seems plausible that a subset $K_0$ of a compact connected full and
locally connected set $K\subset\C$ is always locally connected whenever
it contains any interior component which it meets.

We can now apply Carath\'eodory's Theorem to any connected component of
the interior of $K$: its closure is locally connected by
Lemma~\ref{LemLocConnInterior} above, so its boundary is locally connected
by Carath\'eodory's Theorem, and we then have ``internal rays'' with
respect to any base point in the interior: since this interior component
must be simply connected, it has a Riemann map to $\disk$ sending the
base point to the origin, and the inverse of radial lines under this
Riemann map will be internal rays. Carath\'eodory's Theorem then says
that all internal rays land, and the landing points depend continuously
on the angle. Since no two internal rays can land at the same point
(because the closure of the interior component must be full), the
landing points of the rays induce a homeomorphism between $\Circle$ and
the boundary of the interior component. 

In the remainder of this section, we will consider a fixed compact
connected and full set $K\subset\C$. For the moment, we do not require
it to be locally connected, but we will later add this hypothesis. 

\begin{lemma}[Rays Landing at Common Point]
\label{LemRaysLandTogether} \lineclear
If two external rays land at a common point, then they separate $K$.
\end{lemma}
\proof
The angles of the rays cut $\Circle$ into two parts. If the rays do
not separate $K$, then all the external rays with angles in one of
the two parts of $\Circle$ must land at $z$. But the set of external
angles corresponding to the same landing point always has measure
zero by the Theorem of Riesz \cite[Theorem~\MiLemRiesz{} or
\MiLemRieszA]{MiIntro}.
\qed

The following result will be important for us at several places. Its
proof can be found, for example, in Ahlfors~\cite[Theorem~3.5]{Ahl2}.
\begin{theorem}[Lindel\"of's Theorem]
\label{ThmLindeloef} \lineclear 
If there is a curve $\gamma:[0,1)\to \C-K$ which converges to a point
$z\in\partial K$, then there is a unique external ray of $K$ which
lands at $z$ and which is homotopic to $\gamma$ in $\C-K$.
\qedd
\end{theorem}

\begin{definition}[Access to Boundary Point]
\label{DefAccess} \lineclear
Let $z\in\partial K$. An {\em access} of\/ $z$ is a choice, for every
Euclidean disk $D_r$ of radius $r$ around $z$, of a connected
component\/ $V_r$ of\/ $D_r-K$, such that $V_r\subset V_s$ whenever
$r<s$. 
\end{definition}
\remark
It is not true that, for $r<s$, we must have $V_r=V_s\cap D_r$: there
might be a connected component of $K-\{z\}$ which separates $V_s\cap
D_r$, and $V_r$ is one of its parts. The following lemma justifies
the term ``access''.

\begin{lemma}[Ray in Access]
\label{LemAccess} \lineclear
For every access of\/ $z$, there is a curve in $\C-K$ landing at\/
$z$ running entirely through the domains $V_r$ in the definition of
the access, and visiting all of them. Two such curves can never
separate $K$. Exactly one of these curves is an external ray.
\end{lemma}
\proof
For positive integers $k$ and $r_k=1/k$, let $V_k:=V_{r_k}$, and let
$z_k$ be arbitrary points in $V_k$. Since the $V_k$ are open and nested,
there are curves $\gamma_k$ in $V_k$ connecting $z_k$ to $z_{k+1}$.
Together, they form a curve $\gamma$ starting at $z_1$, remaining in
$V_1$ and necessarily converging to $z$, i.e., landing at $z$. This
curve obviously satisfies the given conditions. If two such curves,
without their landing point $z$, were to separate $K$, then $K$ would
have to be disconnected, a contradiction. If the curves, together
with their landing point, surround some part of $K$, then there is a
radius $r>0$ such that both remaining parts of $K$ contain points at
distance greater than $r$ from $z$. But then these curves cannot stay
forever in the same region $V_r$, so they correspond to different
accesses.  Finally, exactly one of these curves is an external ray by
Lindel\"of's Theorem~\ref{ThmLindeloef}.
\qed

\begin{lemma}[Landing of Rays and Disconnecting Points]
\label{LemDisconnect} \lineclear
The number of external rays landing at any point $z\in\partial K$ equals
the number of connected components of\/ $K-\{z\}$. Between any pair of
external rays, there is a connected component, and conversely. 
\qedd
\end{lemma}
\proof
By Carath\'eodory's Theorem~\ref{ThmCaratheodory}, every boundary point
of $K$ is the landing point of at least one external ray. Any pair of
rays landing at $z$ separates $K$ by Lemma~\ref{LemRaysLandTogether}, so
we only have to show the converse. Let $K_1$ and $K_2$ be two components
of $K-\{z\}$ and let $r>0$ be such that both components contain points
at distance $r$ from $z$. For positive integers $k$, let $D_k$ be the
Euclidean disk of radius $1/k$ around $z$. For $k_0>1/r$, $D_{k_0}-K$ is
disconnected because the sets $K_i$ separate it. For any connected
component of $D_{k_0}-K$ containing $z$ on its closure, there is an
access to $z$ (and possibly many), and there must be two accesses to $z$
separating $K_1$ and $K_2$. 

It follows that any finite collection of rays landing at $z$ produces
equally many connected components of $K-\{z\}$ between them, and any
finite collection of connected components gives rise to equally many
rays separating them. The numbers of rays and connected components
are thus either both finite and equal, or they are both infinite. 
\qed

\remark
Even if $K$ is locally connected, the number of external rays landing
at any single point is not necessarily finite or even countable; see
Lemma~\ref{LemBranchCount}. However, in our applications, this number
will always be finite (although not necessarily bounded over the
branch points of $K$): for connected Julia sets of unicritical
polynomials, three or more rays landing at the same point are always
preperiodic or periodic by Thurston's No Wandering Triangles
Theorem~\ref{ThmThurstonTriangles} (except for the rays on the inverse
orbit of the critical point, but if they are not eventually periodic,
then their number can be at most twice the degree); for the Mandelbrot
set, three rays can land together only at Misiurewicz points
\cite[\CorThreeRaysOneFiber]{FiberMandel}, and the number of rays
landing there is always finite. 

From now on, we assume the set $K$ to be locally connected, in addition
to the requirements that it be compact, connected and full. In the
following, we will collect several properties of such sets. 

\begin{definition}[Branch Point]
\label{DefBranch} \lineclear
A {\em branch point} of a compact connected locally connected full set\/
$K$ is a point which is the landing point of at least three external
rays; equivalently, it is a point which disconnects $K$ into at least
three parts.
\end{definition}

\begin{lemma}[Projection onto Interior Components]
\label{LemProjectComp} \lineclear
Let $K_0$ be the closure of a connected component of the interior of\/
$K$ (a compact connected locally connected full subset of $\C$, let $z$
be a point in the interior of\/ $K_0$ and let $z'$ be a point in
$K-K_0$. Then there is a unique point in $\partial K_0$ through which
every curve in $K$ connecting $z'$ to $z$ must run. This point
disconnects $K$ so that $z$ and $z'$ are in different connected
components.
\end{lemma}
\proof
Let $\gamma_1$ and $\gamma_2$ be two curves connecting $z'$ to $z$. Such
a curve meets $\partial K_0$ in a compact set, so starting from $z'$,
there will be a first point when the curve reaches $\partial K_0$.
Replacing the rest of the curve by a curve within the interior of $K_0$
landing at the same point (which is possible by Carath\'eodory's
Theorem~\ref{ThmCaratheodory} since $\partial K_0$ is locally connected
by Lemma~\ref{LemLocConnInterior}), we may assume that the curve meets
$\partial K_0$ once. If the curves $\gamma_i$ meet $\partial K_0$ in
different points, then the curves, together with $K_0$, enclose a subset
of $\C$ containing boundary points of $K_0$ in its interior. Since
$\partial K_0\subset\partial K$, this contradicts the assumption that
$K$ is full.

Denoting the unique boundary point thus constructed by $\tilde z$, we
now claim that $K-\tilde z$ is disconnected. If it is not, we prove
that it is still arcwise connected: $K-\tilde z$ is locally arcwise
connected, so the set of points in the path component of $z$ is open.
Any limit point different from $\tilde z$ in the path component is
also within the path component because the limit point has a path
connected neighborhood in $K$. Therefore, we can connect $z$ to $z'$
by a path within $K-\{\tilde z\}$, contradicting uniqueness of
$\tilde z$.
\qed

\remark
This way, we obtain a canonical {\em projection} (which is a
retraction) of $K$ onto $K_0$: this projection is the identity on
$K_0$, and outside of $K_0$ it maps to $\partial K_0$ by the
construction above. It is not hard to see that this projection is
continuous. It is locally constant on $K-K_0$. This projection has been
introduced by Douady and Hubbard in \cite[Expos\'e~II.5] {Orsay}. 
From Lemma~\ref{LemDisconnect}, it follows that the projection image of
any point $z'\in K-K_0$ is the landing point of at least two external
rays.

\begin{corollary}[Projection Images Countable]
\label{CorProjectCount} \lineclear
Let $K_0$ be the closure of a connected component of the interior
of\/ $K$ (compact connected locally connected and full). The projection
of\/ $K-K_0$ onto $K_0$ takes images in at most countably many points.
\end{corollary}
\proof
Every image point is the landing point of at least two external rays,
so it separates an open set of external angles from $K_0$. Different
projection points obviously separate different sets of external
angles. The total sum of external angles thus separated is finite, so
the number of projection points must be countable.
\qed

\begin{lemma}[Countably Many Branches]
\label{LemBranchCount} \lineclear
The number of connected components of\/ $K-\{z\}$ is always countable
when $K$ is compact connected locally connected and full. However, the
external angles of $z$ may form a Cantor set.
\end{lemma}
\proof
For any connected component of $K-\{z\}$, pick a point within and let
$\alpha$ be an external angle of this point. Then, by continuity of
landing points, external rays at angles sufficiently close to $\alpha$
land in the same connected component of $K-\{z\}$, so every connected
component takes up an open set of external angles in $\Circle$.

It follows that the external angles of $z$ form the complement of a
dense open subset of $\C$, and if the connected components of
$K-\{z\}$ are arranged so that between any two of them there is
always another one, then the open subsets corresponding to any
connected component always have disjoint closures, and their
complement is a Cantor set.
\qed
\remark
As mentioned above, in all the cases of interest to us the number of
rays landing at a single point will be finite.

\begin{lemma}[Countably Many Branch Points]
\label{LemCountableBranchPts} \lineclear
The number of branch points of and compact connected locally connected
and full subset of\/ $\C$ is countable.
\end{lemma}
\proof
This is a special case of Lemma~\ref{LemBranchCountable}, the proof
of which was self-contained.
\qed

\bigskip
\noindent Dierk Schleicher\\
\medskip
%\parbox[t]{65mm}{ 
Fakult\"at f\"ur Mathematik\\ Technische
Universit\"at\\ Barer Stra{\ss}e 23\\ D-80290 M\"unchen, Germany\\
{\sl dierk$@$mathematik.tu-muenchen.de}
% }
%\hfill and \hfill
%\parbox[t]{65mm}{ Institute for Mathematical Sciences \\ State
%University of New York \\ Stony Brook, NY 11794-3660\\ {\sl
%dierk$@$math.sunysb.edu} }

%\vfill
%\hfill Revision of: 11. September 1998

\end{document}

%% file: LaTeX-top.tex
\newfont{\bbf}{msbm10 at 12pt}

\def\Q{\mbox{\bbf Q}}
\def\Z{\mbox{\bbf Z}}
\def\R{\mbox{\bbf R}}
\def\C{\mbox{\bbf C}}
\def\Cbar{\overline{\cz}}

\def\disk{\mbox{\,\bbf D}}

\def\diskbar{\ovl{\disk}}

\def\Circle{\mbox{\bbf S}^1}

\def\Sphere{\mbox{\bbf S}^2}

\def\cz{\mbox{\bbf C}}

\def\ovl{\overline}
\def\phi1{\phi}
\def\phi{\varphi}
\def\eps{\varepsilon}

\def\*{ {\mbox{\tt $\star$}} }
\def\0{ {\mbox{\tt 0}} }
\def\1{ {\mbox{\tt 1}} }
\def\2{ {\mbox{\tt 2}} }
\def\3{ {\mbox{\tt 3}} }
\def\4{ {\mbox{\tt 4}} }

\def\St#1#2 {{\small\tt\rule{0pt}{0pt}_{\mbox{#1}}^{\mbox{#2}} }}
\def\<{\prec}
\def\>{\succ}

\ifx\ThmNum\undefined 
   \newtheorem{theorem}{Theorem}[section]
   \def\newsection#1{\setcounter{theorem}{0} \section{#1}}
\else 
   \newtheorem{theorem}{Theorem} 
   \def\newsection#1{\section{#1}}
\fi

\newtheorem{proposition}[theorem]{Proposition}
\newtheorem{lemma}[theorem]{Lemma}
\newtheorem{definition}[theorem]{Definition}

\newtheorem{corollary}[theorem]{Corollary}

\def\proof{\par\medskip\noindent {\sc Proof. }}
\def\proofof #1 {\par\medskip\noindent {\sc Proof of #1. }}

\def\sketchof #1 {\par\medskip\noindent {\sc Sketch of proof of #1. }}
\def\Box{\framebox[10pt]{\rule{0pt}{3pt}}}
\def\qed{\hfill $\Box$ \medskip \par}
\def\qedd{\rule{0pt}{0pt}\hfill $\Box$}
\def\remark{\par\medskip \noindent {\sc Remark. }}
\def\lineclear{\rule{0pt}{0pt}\par\noindent}

\def\LabelCaption#1#2{
   \centerline{\parbox{\captionwidth}{   %\vspace{-10mm} 
   \caption{\sl #2}  \label{#1} } }
}

\def\reminder #1 {{\sf #1}}
\def\hide #1 {}
\long\def\longhide #1 {}

%% file: imsmark.tex
\def\IMSmarkvadjust{0 pt}
\def\IMSmarkhadjust{0 pt}
\def\SBIMSMark#1#2#3{
 \font\SBF=cmss10 at 10 true pt
 \font\SBI=cmssi10 at 10 true pt
 \setbox0=\hbox{\SBF Stony Brook IMS Preprint \##1}
 \setbox2=\hbox to \wd0{\hfil \SBI #2}
 \setbox4=\hbox to \wd0{\hfil \SBI #3}
 \setbox6=\hbox to \wd0{\hss
             \vbox{\hsize=\wd0 \parskip=0pt \baselineskip=10 true pt
                   \copy0 \break%
                   \copy2 \break% 
                   \copy4 \break}}
 \dimen0=\ht6   \advance\dimen0 by \vsize \advance\dimen0 by 8 true pt
                \advance\dimen0 by -\pagetotal
	        \advance\dimen0 by \IMSmarkvadjust
 \dimen2=\hsize \advance\dimen2 by .25 true in
	        \advance\dimen2 by \IMSmarkhadjust

%
%   Check for publication info
%
%  \newread\jref
  \openin2=publishd.tex
  \ifeof2\setbox0=\hbox to 0pt{}
  \else 
     \setbox0=\hbox to 3.1 true in{
                \vbox to \ht6{\hsize=3 true in \parskip=0pt  \noindent  
                {\SBI Published in modified form:}\hfil\break
                \input publishd.tex 
                \vfill}}
  \fi
  \closein2
  \ht0=0pt \dp0=0pt
 \ht6=0pt \dp6=0pt
 \setbox8=\vbox to \dimen0{\vfill \hbox to \dimen2{\copy0 \hss \copy6}}
 \ht8=0pt \dp8=0pt \wd8=0pt
 \copy8
 \message{*** Stony Brook IMS Preprint #1, #2. #3 ***}
}